\documentclass[11pt]{article}
\usepackage[utf8]{inputenc}
\usepackage{amsmath,amssymb,amsthm,amsfonts}
\usepackage{enumitem}
\usepackage{needspace}
\setlist[enumerate]{label=(\alph*),font=\normalfont}
\usepackage{geometry}
\usepackage{authblk}
\usepackage{tikz}
\usepackage[hidelinks]{hyperref}
\numberwithin{equation}{section}
\theoremstyle{plain}
\newtheorem{theorem}{Theorem}[section]
\newtheorem{proposition}[theorem]{Proposition}
\newtheorem{lemma}[theorem]{Lemma}

\theoremstyle{definition}

\theoremstyle{remark}
\newtheorem{remark}[theorem]{Remark}
\newtheorem*{unnumberedremark}{Remark}
\newtheorem{example}[theorem]{Example}

\title{Isomorphisms between Leavitt algebras and their matrix rings with degree bounds on generating matrices}

\author{Huynh Viet Khanh}
\author{Nguyen Minh Thanh}
\affil{Department of Mathematics and Informatics, HCMC University of Education\\280 An Duong Vuong Str., Ho Chi Minh City, Vietnam\\\textbf{\texttt{khanhhv@hcmue.edu.vn}}\quad\textbf{\texttt{minhthanhn228@gmail.com}}}
\date{}

\begin{document}
\maketitle

\begin{abstract}
We answer a question of Abrams, \'Anh and Pardo on degree bounds for generators of matrix rings over Leavitt algebras. For a field $K$, an integer $n\ge2$ and a positive integer $d$ relatively prime to $n-1$, we construct $n$ matrices that, together with their adjoints, satisfy the Leavitt relations and generate $M_d(L_K(1,n))$. In the grading in which real edges have degree $-1$ and ghost edges have degree $1$, the homogeneous components of entries of the matrices have degree at most $2$, and those of their adjoints have degree at least $-2$. We realize $L_K(1,n)$ as the Leavitt path algebra of the rose $R_n$, use a partition of the row indices to assign leaves of full $n$-ary trees to matrix positions, and choose row and column depths to control the degrees after expansion.
\end{abstract}



\section{Introduction}\label{sec:introduction}

For a field $K$ and an integer $n\ge2$, let $R=L_K(1,n)$ be the unital $K$-algebra generated by $x_1,\ldots,x_n,y_1,\ldots,y_n$, subject to
$$
x_i y_j=\delta_{ij}1\quad\text{and}\quad\sum_{i=1}^{n}y_i x_i=1.
$$

In \cite[Th.~4.14]{AAP}, Abrams, \'Anh and Pardo proved that $R\cong M_d(R)$ for $\gcd(d,n-1)=1$. For $d\ge2$, their $n$ matrices and their duals have entry degrees $\le d-1$ and $\ge1-d$, respectively, in the grading $\deg x_i=1$, $\deg y_i=-1$. They ask in \cite[Sec.~6, p.~131]{AAP} whether $2$ and $-2$ suffice for a generating set of size $2n$.

\begin{figure}[!ht]
\centering
\begin{tikzpicture}[
    x=1.50pt, y=1.50pt,
    line width=.56pt, line cap=round, line join=round,
    omitted loop/.style={dash pattern=on .48pt off 2.1pt, line cap=butt}
]
\draw[omitted loop] (-5.5,.6) .. controls (-26.4,4.9) and (-3,28.3) .. (1.2,7.3);
\draw[omitted loop] (-5,-5.5) .. controls (-21,-14.7) and (-21.2,14.2) .. (-5.5,5.8);
\draw[omitted loop] (1.6,-5.5) .. controls (-1.7,-29.6) and (-26.9,-5) .. (-5.5,-.7);
\draw[omitted loop] (6.5,-5.5) .. controls (18.6,-23.5) and (-14.9,-23.5) .. (-2.7,-5.5);

\draw (9.2,5.8) .. controls (24.9,14.2) and (24.8,-14.7) .. (8.7,-5.5);
\draw (2.5,7.3) .. controls (6.8,28.3) and (30.1,4.9) .. (9.2,.6);
\draw (-3.8,7.3) .. controls (-12.2,22.9) and (16,22.9) .. (7.6,7.3);
\draw (9.2,-.7) .. controls (30.7,-5) and (5.4,-29.6) .. (2.2,-5.5);

\draw (12.3,-9.3) -- (8.7,-5.5) -- (13.8,-6.6);
\draw (13.9,3.1) -- (9.2,.6) -- (14.4,.1);
\draw (8.7,12.4) -- (7.6,7.3) -- (11.4,10.9);
\draw (1.3,-10.6) -- (2.2,-5.5) -- (4.4,-10.2);

\fill (0,0) circle[radius=2.8pt];
\node[inner sep=0pt] at (4.5,3.1) {$v_0$};
\node[inner sep=0pt] at (28.4,0) {$e_1$};
\node[inner sep=0pt] at (23,21.8) {$e_2$};
\node[inner sep=0pt] at (.3,26.9) {$e_3$};
\node[inner sep=0pt] at (23.3,-21.7) {$e_n$};
\end{tikzpicture}
\setlength{\abovecaptionskip}{6pt}
\caption{The rose $R_n$.}
\label{fig:rose}
\end{figure}
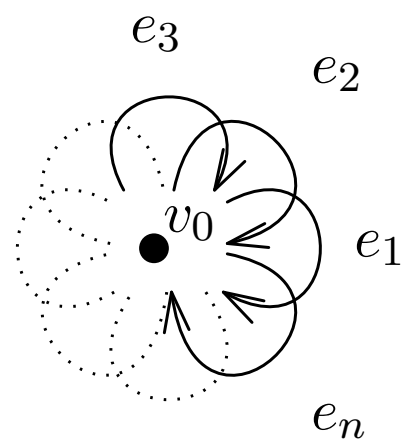

Following Abrams, Ara and Siles Molina \cite[Notation~1.3.1]{AAS2017}, let $R_n$ be the rose with $n$ petals, consisting of a vertex $v_0$ and loops $e_1,\ldots,e_n$ (Figure~\ref{fig:rose}).

It follows from \cite[Prop.~1.3.2]{AAS2017} that $L_K(1,n)\cong L_K(R_n)$. With the generator names in \cite{AAP}, this identification is given by $1\mapsto v_0$, $x_i\mapsto e_i^*$ and $y_i\mapsto e_i$, where $e_i^*$ is the ghost edge of the real edge $e_i$. In $L_K(R_n)$ the Cuntz--Krieger relations (CK1) and (CK2) read
$$
e_i^*e_j=\delta_{ij}v_0\quad\text{and}\quad\sum_{i=1}^{n}e_i e_i^*=v_0.
$$

We regard $v_0$ as the path of length zero. The $K$-linear involution $*$ fixes the scalars and $v_0$, interchanges $e_i$ and $e_i^*$, and reverses products. For a real path $\lambda=e_{i_1}\cdots e_{i_t}$ of length $t\ge1$, the corresponding ghost path is $\lambda^*=e_{i_t}^*\cdots e_{i_1}^*$. For a matrix $C=(c_{pq})$, its adjoint has entries $(C^*)_{pq}=(c_{qp})^*$.

We use the $\mathbb Z$-grading in which $\deg e_i=-1$ and $\deg e_i^*=1$. Matrix units have degree zero. For paths $\alpha,\beta$ in $R_n$, the monomial $\alpha\beta^*$ has degree $|\beta|-|\alpha|$, where $|\lambda|$ denotes the length of a path $\lambda$. This is the negative of its degree under the standard grading of $L_K(R_n)$ described in \cite[Cor.~2.1.5(iii)]{AAS2017}. If $R_t$ denotes the homogeneous component of degree $t$ in our grading, set
$$
R_{\le2}=\mathop{\textstyle\bigoplus}\limits_{t\le2}R_t\quad\text{and}\quad R_{\ge-2}=\mathop{\textstyle\bigoplus}\limits_{t\ge-2}R_t.
$$

The following theorem answers the question with the dual matrices realized as adjoints.

\begin{theorem}\label{thm:main}
Let $K$ be a field, $n\ge2$, and $d\ge1$, with $\gcd(d,n-1)=1$. There exist matrices $X_1,\ldots,X_n$ in $M_d(R_{\le2})$ such that, with $Y_i=X_i^*$,
$$
X_iY_j=\delta_{ij}I_d\quad\text{and}\quad\sum_{i=1}^{n}Y_iX_i=I_d,
$$
and the matrices $X_1,\ldots,X_n,Y_1,\ldots,Y_n$ generate $M_d(R)$ as a $K$-algebra. In particular, $Y_i\in M_d(R_{\ge-2})$.
\end{theorem}

The bounds apply to each homogeneous component of a matrix entry. For a path monomial $\alpha\beta^*$, the degree is the signed difference $|\beta|-|\alpha|$, and the total path length $|\alpha|+|\beta|$ may exceed $2$. By Remark~\ref{rem:distinct-generators}, the resulting generating set has exactly $2n$ elements.

A tuple $(X_1,\ldots,X_n)$ in $M_d(R)$ is called a \emph{Leavitt generating family} if the matrices $X_i$ and $X_i^*$ satisfy the two displayed relations of Theorem~\ref{thm:main} and together generate $M_d(R)$. The assignments $e_i^*\mapsto X_i$ and $e_i\mapsto X_i^*$ induce a unital $K$-linear $*$-homomorphism $R\to M_d(R)$. The generation assumption implies surjectivity. Injectivity follows from the simplicity of $R$ recorded in \cite[Prop.~1.1(3)]{AAP}. Thus the assignments determine an isomorphism.

Put $N=n-1$. In the matrices of Abrams, \'Anh and Pardo, an entry $1$ at $(p,q)$ defines the directed edge $p\to q$ on the row indices. The vertex sets $S_1,S_2$ of the two resulting paths form the partition used in Section~\ref{sec:partition} to count monomial positions. Their matrices generate for every assignment of the fixed monomial list respecting this partition. For $2\le a<n$ and $\gcd(a,N)=1$, the matrices in Proposition~\ref{prop:canonical} generate $M_a(R)$ for every color-preserving assignment of leaf adjoints of a compatible full $n$-ary tree.

In Section~\ref{sec:transport}, a nonzero entry obtained from $c_{pq}\in R_t$ has degree $t+d_p-d_q$, where $d_p$ and $d_q$ are the row and column depths. Lemma~\ref{lem:lifting} uses a directed graph of possible degree-two entries to select indices for further expansion. Under its support condition, the bound $2$ transfers from size $e$ to every size $d\ge1$ with $d\equiv e\pmod N$.

Section~\ref{sec:direct} verifies the degree bound and support condition for the direct and composed constructions. Section~\ref{sec:remaining} chooses the composition sizes for prime $n$. For composite $n$, the graded isomorphism $M_d(R)\cong M_{dp^\ell}(R)$ with $p$ a prime divisor of $n$ reduces the proof to $dp^\ell\equiv1\pmod N$ or to a class covered by the direct construction.

\section{The partition and generating matrices}\label{sec:partition}\label{sec:canonical}

Assume $n\ge3$, and let $2\le a<n$, with $\gcd(a,N)=1$. For $x\in\mathbb R$, let $\lfloor x\rfloor$ denote the greatest integer $\le x$. Write
$$
n=Qa+r,\qquad Q=\Bigl\lfloor\frac{N-1}{a}\Bigr\rfloor,\qquad \delta=r-1\quad\text{and}\quad s=a-\delta.
$$
Since $2\le a\le N$ and $\gcd(a,N)=1$, we have $a\le N-1$ and $Q\ge1$. The remainder $r-2=N-1-Qa$ satisfies $0\le r-2<a$. Also $\delta=N-Qa$, so $\gcd(\delta,a)=1$ and $\delta\ne a$. Hence $2\le r\le a$. For a unit $h$ modulo $N$, define $f(h)\in\{1,\ldots,N-1\}$ by $hf(h)\equiv-1\pmod N$. Then $f(f(h))=h$ for $1\le h<N$ and $\gcd(h,N)=1$.

Choose $m\in\{1,\ldots,a-1\}$ with $m\delta\equiv1\pmod a$. Let $\sigma_i$ be the representative of $1+(i-1)s$ modulo $a$ in $\{1,\ldots,a\}$, and set
$$
S_1=\{\sigma_1,\ldots,\sigma_m\}\quad\text{and}\quad S_2=\{\sigma_{m+1},\ldots,\sigma_a\}.
$$
For $1\le u\le n$, put $\widehat u=1+((u-1)\bmod a)$ and define
$$
B_a=\{u\in\{1,\ldots,n\}:\widehat u\in S_1\}\quad\text{and}\quad \beta_a=|B_a|.
$$
For fixed $a$, write $B=B_a$ and $\beta=\beta_a$.

We use $S_1,S_2$ and $B_a$ for the sets arising from the partition construction in \cite[Defs.~1.4, 1.6 and 4.1]{AAP}.

\begin{lemma}\label{lem:partition}
The sequence $(\sigma_1,\ldots,\sigma_a)$ is a permutation of $\{1,\ldots,a\}$, with $\sigma_1=1$, $\sigma_m=\delta$, $\sigma_{m+1}=a$ and $\sigma_a=r$. Moreover, $af(a)+1=mN$, $\beta=f(a)+1$ and $Q+1\le\beta\le n-Q-1$.
\end{lemma}

\begin{proof}
Since $s\equiv-\delta\pmod a$ and $\gcd(\delta,a)=1$, the $\sigma_i$ form a permutation. We have $\sigma_1=1$ and
$$
\sigma_m\equiv\delta,\qquad \sigma_{m+1}\equiv0,\qquad \sigma_a\equiv1+\delta=r\pmod a.
$$
The representatives in $\{1,\ldots,a\}$ are therefore the asserted endpoints.

Put $f_0=|S_1\cap\{1,\ldots,r\}|$. For a condition $P$, its indicator $\mathbf1_{\{P\}}$ is $1$ if $P$ holds and $0$ otherwise. The representatives $\sigma_i\in\{1,\ldots,a\}$ satisfy
$$
\sigma_{i+1}-\sigma_i=-\delta+a\mathbf1_{\{\sigma_i\le\delta\}}\qquad\text{for }1\le i\le m.
$$
Since $r=\delta+1\in S_2$, we have $\sum_{i=1}^m\mathbf1_{\{\sigma_i\le\delta\}}=f_0$. Summing the successive differences yields
$$
a-1=\sigma_{m+1}-\sigma_1=-m\delta+af_0,\quad\text{so}\quad a(f_0-1)=m\delta-1.
$$
Counting the elements of $B$ in the $Q$ full blocks of labels and the last block of length $r$, we obtain $\beta=Qm+f_0$. Since $N=Qa+\delta$,
$$
a(\beta-1)=aQm+a(f_0-1)=mN-1.
$$
Since $1\le m<a$, we have $0<\beta-1<N$. Thus $\beta-1=f(a)$ and $af(a)+1=mN$.

The $Q+1$ indices $1+ja$, with $0\le j\le Q$, have reduced row $1$ and belong to $B$. The $Q+1$ distinct indices $a,2a,\ldots,Qa,n$ have reduced rows in $\{a,r\}\subseteq S_2$ and belong to $B^c$. Therefore $Q+1\le\beta\le n-Q-1$.
\end{proof}

\begin{unnumberedremark}[Christoffel words]
Put $q=a-m$. Since $m\delta\equiv1\pmod a$, we have $\gcd(m,q)=1$. The lower Christoffel word with $m$ zeros and $q$ ones is $\epsilon_1\cdots\epsilon_a$, where
$$
\epsilon_j=\Bigl\lfloor\frac{jq}{a}\Bigr\rfloor-\Bigl\lfloor\frac{(j-1)q}{a}\Bigr\rfloor\qquad\text{for }1\le j\le a.
$$
Let $\rho_j$ be the least nonnegative residue of $(j-1)q$ modulo $a$. The equivalences
$$
\epsilon_j=0\iff\rho_j+q<a\iff\rho_j<m
$$
express the residue criterion in \cite[Sec.~4.1, Eq.~(4)]{DKKS2024} for the letter $0$. Since $sq\equiv1\pmod a$, we have $\rho_{\sigma_i}=i-1$ for $1\le i\le a$. Hence $S_1$ is the set of zero positions.
\end{unnumberedremark}

A scalar $1$ at $(p,q)$ in a matrix of Proposition~\ref{prop:canonical} determines $p\to q$ with $q\equiv p+s\pmod a$. Rows $\delta$ and $r$ have no scalar entries, so $\delta\to a$ and $r\to1$ are absent. By Lemma~\ref{lem:partition}, the scalar edges form the paths
$$
(\sigma_1,\ldots,\sigma_m)\quad\text{and}\quad (\sigma_{m+1},\ldots,\sigma_a).
$$
The initial vertices are $\sigma_1=1$ and $\sigma_{m+1}=a$. Columns $1$ and $a$ contain the ghost edges $e_1^*,\ldots,e_n^*$ and the adjoints of the assigned leaves, respectively.

The \emph{monomial positions} $(i,p)$ index the entries $(X_i)_{pa}$ with $i=Q+2$ and $r-1\le p\le a$, or $Q+3\le i\le n$ and $1\le p\le a$. Such a position has color $c\in\{1,2\}$ if $p\in S_c$. For $1\le p\le a$, the number in row $p$ is
$$
C_a(p)=n-Q-2+\mathbf1_{\{p\ge r-1\}}.
$$

A finite full $n$-ary tree is a finite prefix-closed set of paths in $R_n$ with root $v_0$ and children $\{\nu e_u:1\le u\le n\}$ at each internal vertex $\nu$. Its leaves form an $n$-ary leaf set, as defined by H.~V.~Khanh and V.~H.~Thanh in \cite[Sec.~7]{KT2026}. The path $\nu e_u$ has color $c$ if $\widehat u\in S_c$. The tree is \emph{compatible} if it has $a-1$ internal vertices and each nonroot internal path consists of edges $e_u$ with $u\in B$.

The count of monomial positions in the first part, obtained in \cite[Lem.~4.7 and Prop.~4.9]{AAP}, is expressed below in terms of $\beta$.

\begin{lemma}\label{lem:counts}
The total numbers of monomial positions in the two parts are
$$
\sum_{p\in S_1}C_a(p)=(a-1)(\beta-1)+1\quad\text{and}\quad \sum_{p\in S_2}C_a(p)=(a-1)(n-\beta).
$$
These are also the numbers of leaves of the respective colors in every compatible tree. Put $\Delta=f(a)-Q\ge0$. At the initial vertices $1$ and $a$ of the two directed paths,
$$
C_a(1)=N-Q-1+\mathbf1_{\{r=2\}}\quad\text{and}\quad C_a(a)=(n-\beta)+\Delta.
$$
\end{lemma}

\begin{proof}
Put $e_0=|S_1\cap\{r-1,\ldots,a\}|$ and $f_0=|S_1\cap\{1,\ldots,r\}|$. These two intervals cover $\{1,\ldots,a\}$ and intersect in $\{r-1,r\}$, of which only $r-1$ lies in $S_1$. Thus $e_0+f_0=m+1$. Using $m(r-1)=a(f_0-1)+1$ and $\beta=Qm+f_0$ from the proof of Lemma~\ref{lem:partition}, we obtain
$$
\begin{aligned}
\sum_{p\in S_1}C_a(p)
&=m(n-Q-2)+e_0\\
&=mQ(a-1)+m(r-1)-f_0+1\\
&=(a-1)(Qm+f_0-1)+1=(a-1)(\beta-1)+1.
\end{aligned}
$$
Subtracting the first count from the total $\sum_{p=1}^a C_a(p)=a(n-Q-2)+(a-r+2)=1+(a-1)N$ yields $\sum_{p\in S_2}C_a(p)=(a-1)(n-\beta)$.

Each internal path $\nu$ has $\beta$ children $\nu e_u$ with $u\in B$ and $n-\beta$ with $u\notin B$. Among all these children, the $a-2$ nonroot internal vertices have color one. The leaf counts are therefore $(a-1)\beta-(a-2)$ and $(a-1)(n-\beta)$, respectively.

By the definition of $C_a$, we have $C_a(1)=N-Q-1+\mathbf1_{\{r=2\}}$ and $C_a(a)=n-Q-1=(n-\beta)+f(a)-Q$. By Lemma~\ref{lem:partition}, $f(a)=\beta-1\ge Q$.
\end{proof}

For $h$ satisfying the same hypotheses, $Q_h,r_h,B_h,\beta_h,C_h$ denote the corresponding quantities with $a$ replaced by $h$.

Write $E_{pq}$ for the standard matrix units of $M_a(R)$. Fix a compatible tree $T$ and use Lemma~\ref{lem:counts} to assign its leaves bijectively to monomial positions of the same color. If $\lambda=\nu e_u$ occupies $(i,j)$, put $a_{i,j}=\lambda^*$, with $j,\widehat u\in S_c$ for some $c\in\{1,2\}$.

We use the nonzero entries in column $1$ and the scalar entries given in \cite[Sec.~3]{AAP}. We first obtain $E_{pp}$ and then $E_{pq}$ for $p,q\in S_c$, as in \cite[Props.~3.10 and 3.14]{AAP}. Induction on internal vertices extends the reduction of powers of $x_1$ used in \cite[Prop.~4.13]{AAP}.

\begin{proposition}\label{prop:canonical}
For any such bijection, the matrices
$$
\begin{aligned}
X_i&=\sum_{j=1}^{a}e_{(i-1)a+j}^*E_{j1},&&\text{for }1\le i\le Q,\\
X_{Q+1}&=\sum_{j=1}^{r}e_{Qa+j}^*E_{j1}+\sum_{j=r+1}^{a}E_{j,j-\delta},\\
X_{Q+2}&=\sum_{j=1}^{r-2}E_{j,j+s}+\sum_{j=r-1}^{a}a_{Q+2,j}E_{ja},\\
X_i&=\sum_{j=1}^{a}a_{i,j}E_{ja},&&\text{for }Q+3\le i\le n,
\end{aligned}
$$
form a Leavitt generating family in $M_a(R)$, with empty sums understood as zero. Consequently, $e_i^*\mapsto X_i$ and $e_i\mapsto X_i^*$ induce a unital $K$-linear $*$-isomorphism $R\longrightarrow M_a(R)$.
\end{proposition}

\begin{proof}
For the leaf set $\Lambda$ of a finite full tree, distinct leaves are incomparable under the prefix relation. The prefix cancellation rule in \cite[Lem.~1.2.12(i)]{AAS2017} then implies $\lambda^*\mu=\delta_{\lambda\mu}v_0$ for $\lambda,\mu\in\Lambda$. For each internal path $\nu$, (CK2) at $v_0$ implies
$$
\nu\nu^*=\nu\Bigl(\sum_{u=1}^n e_ue_u^*\Bigr)\nu^*
=\sum_{u=1}^n(\nu e_u)(\nu e_u)^*.
$$
Expanding from $v_0$ to the leaves proves $\sum_{\lambda\in\Lambda}\lambda\lambda^*=v_0=1$.

Let $c_i(p)$ be the column of the unique nonzero entry $\lambda_{i,p}^*$ in row $p$ of $X_i$. For each column $q$, the paths $\lambda_{i,p}$ with $c_i(p)=q$ run once through a full tree's leaves. These are the $e_u$ for $q=1$, just $v_0$ for $2\le q<a$, and the leaves of $T$ for $q=a$. Consequently,
$$
\begin{aligned}
(X_iX_j^*)_{pq}
&=\delta_{c_i(p),c_j(q)}\lambda_{i,p}^*\lambda_{j,q}
=\delta_{ij}\delta_{pq}1,\\
\sum_{i=1}^n X_i^*X_i
&=\sum_{q=1}^a\biggl(\sum_{i,p:\,c_i(p)=q}\lambda_{i,p}\lambda_{i,p}^*\biggr)E_{qq}
=I_a.
\end{aligned}
$$
Put $E_t=\sum_{j=1}^t E_{jj}$, with $E_0=0$. The two ranges of generators $1\le i\le Q+1$ and $Q+2\le i\le n$ use disjoint column ranges $1,\ldots,s$ and $s+1,\ldots,a$. Thus their respective sums of $X_i^*X_i$ are $E_s$ and $I_a-E_s$.

Let $A$ be the $K$-algebra generated by the $X_i$ and $X_i^*$. In particular $I_a,E_s\in A$, and $A$ is closed under $*$. For $1\le i\le a$, set $u_i=[is]_a$ with representatives in $\{1,\ldots,a\}$. Since $\gcd(s,a)=1$, the $u_i$ form a permutation, with $u_1=s$, $u_{a-1}=\delta$ and $u_a=a$. For $1\le i\le a-2$, we have
$$
E_{u_{i+1}}=
\begin{cases}
E_s+X_{Q+2}^* E_{u_i}X_{Q+2},&u_i<\delta,\\
E_s-X_{Q+1}^*(I_a-E_{u_i})X_{Q+1},&u_i>\delta.
\end{cases}
$$
Indeed, when $u_i<\delta$, the added term is $\sum_{j=1}^{u_i}E_{j+s,j+s}$ and $u_{i+1}=u_i+s$. When $u_i>\delta$, the subtracted term is $\sum_{j=u_i+1}^a E_{j-\delta,j-\delta}$ and $u_{i+1}=u_i-\delta$. Induction shows that $E_t\in A$ for all $0\le t\le a$, and hence $E_{pp}=E_p-E_{p-1}\in A$.

At a scalar position $(p,q)$ of $X_i$, we have $E_{pp}X_iE_{qq}=E_{pq}\in A$. The scalar positions determine the two directed paths on $S_1$ and $S_2$ described after Lemma~\ref{lem:partition}. Products along a path and their adjoints yield $E_{pq}\in A$ for $p,q$ in the same part, with $p=q$ already covered. Likewise, for $1\le u\le n$ and $i=\lfloor(u-1)/a\rfloor+1$,
$$
E_{\widehat u,\widehat u}X_iE_{11}=e_u^*E_{\widehat u,1}\in A
\quad\text{and}\quad e_uE_{1,\widehat u}\in A.
$$

We prove $\nu^*E_{1a}\in A$ for each internal path $\nu$, starting at greatest depth. If $\nu e_u$ is a leaf assigned to row $j$ of $X_i$, then $j$ and $\widehat u$ belong to the same part. Therefore
$$
(e_uE_{1,\widehat u})E_{\widehat u,j}(E_{jj}X_iE_{aa})
=e_u(\nu e_u)^*E_{1a}\in A.
$$
If $\nu e_u$ is internal, then $\widehat u\in S_1$ by compatibility and $(\nu e_u)^*E_{1a}\in A$ by induction. In this case,
$$
(e_uE_{1,\widehat u})E_{\widehat u,1}\bigl((\nu e_u)^*E_{1a}\bigr)
=e_u(\nu e_u)^*E_{1a}\in A.
$$
At maximal internal depth only the leaf case occurs. Summing the contributions of the children and using $(\nu e_u)^*=e_u^*\nu^*$ yields
$$
\sum_{u=1}^n e_u(\nu e_u)^*E_{1a}
=\Bigl(\sum_{u=1}^n e_ue_u^*\Bigr)\nu^*E_{1a}
=\nu^*E_{1a}\in A.
$$

At the root $\nu=v_0$, so $E_{1a}\in A$. For $p\in S_1$ and $q\in S_2$, we have $E_{pq}=E_{p1}E_{1a}E_{aq}\in A$, from which it follows by taking adjoints that $E_{qp}\in A$. With all matrix units in $A$, we obtain
$$
E_{p,\widehat u}(e_u^*E_{\widehat u,1})E_{1q}=e_u^*E_{pq}\in A
\qquad\text{for all }u,p,q.
$$
By closure under the adjoint, $e_uE_{pq}\in A$ for all $u,p,q$. These elements and the matrix units generate $M_a(R)$. Since $R$ is simple by \cite[Prop.~1.1(3)]{AAP}, the induced homomorphism is injective and hence an isomorphism.
\end{proof}

Each generator row has exactly one nonzero entry, a path adjoint. For two matrices $A,B$ of this form with $A_{pq}=\alpha^*$ and $B_{qk}=\beta^*$, row $p$ of $AB$ is supported in column $k$, with
$$
(AB)_{pk}=\alpha^*\beta^*=(\beta\alpha)^*.
$$
Since $(\beta\alpha)^*(\beta\alpha)=1$, this entry is nonzero. Substituting one family of Proposition~\ref{prop:canonical} into another replaces $\lambda^*$ by a product of generators, with $I_b$ for $\lambda=v_0$. The row property is preserved under substitution with the lexicographic identification $M_a(M_b(R))\cong M_{ab}(R)$ and under permutations of the generator indices.

\section{Matrix expansions and degree bounds}\label{sec:transport}

For $n\ge2$, Lemma~\ref{lem:forest} expands each index $p$ to $\{(p,\alpha):|\alpha|=d_p\}$. An entry $c_{pq}\in R_t$ expands to entries in $R_{t+d_p-d_q}$. For the matrices of Proposition~\ref{prop:canonical}, these entries follow by prefix comparison in $R_n$. The support condition in Lemma~\ref{lem:lifting} lets us transfer the bound $2$ to every positive dimension in the same residue class modulo $N$.

The isomorphism in Lemma~\ref{lem:forest} iterates the map given by Pardo in \cite[Rem.~2.2]{Pardo2011} after permutations of the matrix coordinates. Isomorphisms associated with ordered leaf sets are described in \cite[Prop.~2.4 and Sec.~7]{KT2026}.

\begin{lemma}\label{lem:forest}
Let $e\ge1$ and $(d_1,\ldots,d_e)\in\mathbb Z_{\ge0}^e$. With $\alpha$ ranging over paths in $R_n$, set
$$
I=\{(p,\alpha):1\le p\le e,\ |\alpha|=d_p\}\quad\text{and}\quad
E=|I|=\sum_{p=1}^e n^{d_p}\equiv e\pmod N.
$$
With rows and columns of $M_E(R)$ indexed by $I$, the map
$$
F:M_e(R)\longrightarrow M_E(R),\qquad
C=(c_{pq})\longmapsto\bigl(\alpha^*c_{pq}\eta\bigr)_{(p,\alpha),(q,\eta)}
$$
is a unital $K$-linear $*$-isomorphism. If $c_{pq}\in R_j$, then $\alpha^*c_{pq}\eta\in R_{j+d_p-d_q}$. For $d_1=\cdots=d_e$, both $F$ and $F^{-1}$ are graded.

Suppose each row of $C$ is zero or has a single nonzero entry $\lambda^*$ for a path $\lambda$. Each row $(r_z)_{z\in I}$ of $F(C)$ then satisfies one of the following.
\begin{enumerate}[leftmargin=*,nosep]
\item Every nonzero $r_z$ is a real path $\xi$ with $|\xi|>0$.
\item $r_{z_0}=\xi^*$ for some $z_0\in I$ and a path $\xi$, with $r_z=0$ for $z\ne z_0$.
\end{enumerate}
\end{lemma}

\begin{proof}
For $t\ge0$, the paths of length $t$ are the leaves of the tree $\{\nu:|\nu|\le t\}$. By prefix cancellation and (CK2), as in Proposition~\ref{prop:canonical},
$$
\alpha^*\eta=\delta_{\alpha\eta}v_0\quad\text{for }|\alpha|=|\eta|=t
\qquad\text{and}\qquad
\sum_{|\alpha|=t}\alpha\alpha^*=v_0.
$$
Put $U=(\delta_{pq}\alpha^*)_{(p,\alpha),q}\in M_{E\times e}(R)$ and $V=U^*$. Then
$$
\begin{aligned}
(UV)_{(p,\alpha),(q,\eta)}
  &=\delta_{pq}\alpha^*\eta=\delta_{pq}\delta_{\alpha\eta}v_0,\\
(VU)_{pq}
  &=\delta_{pq}\sum_{|\alpha|=d_p}\alpha\alpha^*=\delta_{pq}v_0.
\end{aligned}
$$
Thus $UV=I_E$ and $VU=I_e$. By the entry formula, $F(C)=UCV$. Since $V=U^*$, this is a unital $K$-linear $*$-isomorphism with inverse $F^{-1}(D)=VDU$. For $c_{pq},D_{(p,\alpha),(q,\eta)}\in R_j$, the grading yields
$$
\alpha^*c_{pq}\eta\in R_{j+d_p-d_q}\quad\text{and}\quad
\alpha D_{(p,\alpha),(q,\eta)}\eta^*\in R_{j+d_q-d_p}.
$$
Both shifts vanish when $d_1=\cdots=d_e$.

Fix $(p,\alpha)\in I$ and write $r_z=F(C)_{(p,\alpha),z}$ for $z\in I$. If $c_{pq}=\lambda^*$ is the sole nonzero entry of row $p$, put $\mu=\lambda\alpha$. Then $r_{(q,\eta)}=\mu^*\eta$, and $r_{(q',\eta)}=0$ for $q'\ne q$. Apply prefix cancellation as in \cite[Lem.~1.2.12(i)]{AAS2017}.
\begin{enumerate}[leftmargin=*,nosep]
\item If row $p$ of $C$ is zero, then $r_z=0$ for all $z\in I$. If row $p$ is nonzero and $|\mu|<d_q$, the $n^{d_q-|\mu|}$ nonzero entries are
$$
r_{(q,\mu\xi)}=\xi\in R_{|\mu|-d_q}\qquad\text{for }|\xi|=d_q-|\mu|>0.
$$
\item If row $p$ is nonzero and $|\mu|\ge d_q$, write $\mu=\eta_0\xi$ with $|\eta_0|=d_q$ and set $z_0=(q,\eta_0)$. Then $r_z=\delta_{z,z_0}\xi^*$ for $z\in I$, with $|\xi|=|\mu|-d_q$ and $\xi=v_0$ if $|\mu|=d_q$.
\end{enumerate}
Thus each expanded row has at most one nonzero entry of degree two, of the form $\xi^*$ with $|\xi|=2$.
\end{proof}

For the degree-two components $A_i^{(2)}$, Lemma~\ref{lem:lifting} assumes $(A_i^{(2)})_{pq}\in\operatorname{span}_K\{\tau^*:|\tau|=2\}$. With the generator names in \cite{AAP}, the corresponding quadratic monomials satisfy $x_ux_v=(e_ve_u)^*$. The row support $\mathcal S_p$ consists of the pairs $(q,\tau)$ with nonzero coefficient of $\tau^*$ in $(A_i^{(2)})_{pq}$ for some $i$.

\Needspace{14\baselineskip}
\begin{lemma}\label{lem:lifting}
Let $e\ge1$, and let $(A_1,\ldots,A_n)$ be a Leavitt generating family in $M_e(R)$ with $A_i\in M_e(R_{\le2})$ for each $i$. For each $i$, let $A_i^{(2)}$ be the matrix of homogeneous components of degree two. Assume the following.
\begin{enumerate}[leftmargin=*,nosep]
\item For $1\le i\le n$ and $1\le p,q\le e$, the entry $(A_i^{(2)})_{pq}$ has the form $\sum_{|\tau|=2}c_{i,pq}^{\tau}\tau^*$, with $c_{i,pq}^{\tau}\in K$, where $\tau$ runs over paths in $R_n$.
\item For $1\le p\le e$, the set
$$
\mathcal S_p=\bigcup_{i=1}^n\{(q,\tau):1\le q\le e,\ |\tau|=2,\ c_{i,pq}^{\tau}\ne0\}
$$
satisfies $|\mathcal S_p|\le n$.
\end{enumerate}
Then, for every $d\ge1$ with $d\equiv e\pmod N$, there is a Leavitt generating family $(B_1,\ldots,B_n)$ in $M_d(R)$ with $B_i\in M_d(R_{\le2})$ for each $i$.
\end{lemma}

\begin{proof}
Uniform expansions at depths $h$ and $j$ take $M_e(R)$ and $M_d(R)$ to $M_{en^h}(R)$ and $M_{dn^j}(R)$, respectively. To match the sizes, we expand $b=(dn^j-en^h)/N$ indices of $M_{en^h}(R)$ by one further level while preserving the degree bound. We then apply the inverse of the uniform expansion from $M_d(R)$.

Fix $d\ge1$ with $d\equiv e\pmod N$. Choose $\ell\in\mathbb Z$ with $en^\ell\le d<en^{\ell+1}$ and put $\lambda=(dn^{-\ell}-e)/N\in[0,e)$. Choose $k\ge1$ so that $j=2k-\ell\ge0$ and, if $\lambda>0$, $\lambda n^k\ge e$. Put $h=2k$. Since $dn^j\equiv en^h\pmod N$, we have
$$
b=\frac{dn^j-en^h}{N}=\lambda n^h\in\mathbb Z,\qquad
\begin{cases}
b=0,&\lambda=0,\\
en^k\le b<en^h,&\lambda>0.
\end{cases}
$$

Let $F_h$ be the uniform depth-$h$ expansion of Lemma~\ref{lem:forest}, with index set $V_h=\{(p,\alpha):1\le p\le e,\ |\alpha|=h\}$. For $\alpha=\alpha'\xi$ with $|\xi|=2$, gradedness and prefix cancellation yield
$$
\begin{aligned}
\bigl(F_h(A_i)^{(2)}\bigr)_{(p,\alpha),(q,\eta)}
&=\sum_{|\tau|=2}c_{i,pq}^{\tau}(\tau\alpha)^*\eta\\
&=\sum_{\substack{|\tau|=2\\\eta=\tau\alpha'}}c_{i,pq}^{\tau}\xi^*.
\end{aligned}
$$
Define the directed graph $G_h$ on $V_h$ by the edges
$$
(p,\alpha)\xrightarrow{(p,q,\tau)}(q,\tau\alpha')
\qquad\text{for }(q,\tau)\in\mathcal S_p.
$$
Thus a nonzero entry of $F_h(A_i)^{(2)}$ at $(z,z')$ implies $z\to z'$ in $G_h$.

Write $\operatorname{pref}_s(\alpha)$ for the prefix of length $s$, with $\operatorname{pref}_0(\alpha)=v_0$. For a walk from $(p_0,\alpha)$ with successive labels $(p_{m-1},p_m,\tau_m)$, write $(p_m,\alpha_m)$ for its vertex after $m$ steps, with $\alpha_0=\alpha$. For $1\le m\le k$ we have $h-2(m-1)\ge2$, so the recursion yields
$$
\alpha_m=\tau_m\operatorname{pref}_{h-2}(\alpha_{m-1})
=\tau_m\cdots\tau_1\operatorname{pref}_{h-2m}(\alpha).
$$
At $m=k$, the endpoint $(p_k,\tau_k\cdots\tau_1)$ depends only on $p_0$ and the choices $(p_m,\tau_m)\in\mathcal S_{p_{m-1}}$ for $1\le m\le k$. There are $e$ choices for $p_0$ and at most $n$ at each step, so the set $C_h$ of endpoints of length-$k$ walks satisfies $|C_h|\le en^k$.

For $W\subseteq V_h$, let $G_h(W)$ be the set of endpoints of edges from $W$. Then $C_h=G_h^k(V_h)$. Since $G_h$ preserves inclusions,
$$
G_h(C_h)=G_h^k(G_h(V_h))\subseteq G_h^k(V_h)=C_h.
$$
Every vertex on a directed cycle ends a walk of length $k$ along that cycle, so the induced graph on $V_h\setminus C_h$ is acyclic. Set $W_0=C_h$. If $W_r\ne V_h$, choose a sink $z_r$ in the graph induced on $V_h\setminus W_r$, so $G_h(\{z_r\})\subseteq W_r$. Put $W_{r+1}=W_r\cup\{z_r\}$. Inductively,
$$
G_h(W_{r+1})=G_h(W_r)\cup G_h(\{z_r\})\subseteq W_r\subseteq W_{r+1}.
$$
As $|W_r|=|C_h|+r$, for $b>0$ the bounds $|C_h|\le en^k\le b<|V_h|$ allow us to take $W=W_{b-|C_h|}$. For $b=0$, take $W=\varnothing$. In both cases $|W|=b$ and $G_h(W)\subseteq W$.

Let $F_\varepsilon$ be the expansion of Lemma~\ref{lem:forest} with depths $\varepsilon_z=\mathbf1_W(z)$. For a nonzero homogeneous component of degree $t\le2$ at $(z,z')$ in $F_h(A_i)$, each nonzero term after expansion has degree
$$
t+\varepsilon_z-\varepsilon_{z'}\le
\begin{cases}
t+1\le2,&t\le1,\\
2,&t=2.
\end{cases}
$$
For $t=2$, the edge $z\to z'$ and the inclusion $G_h(W)\subseteq W$ imply $\varepsilon_z\le\varepsilon_{z'}$. Thus $F_\varepsilon(F_h(A_i))$ has entries in $R_{\le2}$, and its size is
$$
\sum_{z\in V_h} n^{\varepsilon_z}=|V_h|+N|W|=en^h+bN=dn^j.
$$
Let $F_j:M_d(R)\to M_{dn^j}(R)$ be the uniform depth-$j$ expansion. Since $F_j^{-1}$ is graded,
$$
B_i=(F_j^{-1}\circ F_\varepsilon\circ F_h)(A_i)\in M_d(R_{\le2}).
$$
The composite is a unital $K$-linear $*$-isomorphism, so the $B_i$ form a Leavitt generating family.
\end{proof}

Let $A_1,\ldots,A_n$ be the expanded matrices in Lemma~\ref{lem:forest} when all input rows have at most one nonzero entry, a path adjoint. Each row of $A_i^{(2)}$ then has at most one nonzero entry $\tau^*$ with $|\tau|=2$. By (CK1), these adjoints are linearly independent. Thus each $A_i$ contributes at most one pair $(q,\tau)$ per row, so $|\mathcal S_p|\le n$, as required in Lemma~\ref{lem:lifting}.

\section{Direct and composed constructions}\label{sec:direct}\label{sec:composition}

For the matrices of Proposition~\ref{prop:canonical}, choose $d_q=d_p-2$ along each row-index edge $p\to q$, with $d_a=d_1+1$. By Lemma~\ref{lem:forest}, a nonzero entry arising from $\lambda^*$ at $(p,q)$ has degree $|\lambda|+d_p-d_q$. The scalar positions and column $1$ meet the bound $2$. In column $a$, the bound holds for $|\lambda|\le2+d_a-d_p$. The two cases use compatible trees of height at most $3$ and $5$, respectively.

\begin{proposition}\label{prop:direct}
Let $n\ge3$, and let $2\le a<n$ with $\gcd(a,N)=1$. Put $w=f(a)$. If $a\le w^2+3w+4$ or $N\le w^3$, then, for every field $K$ and every positive integer $d\equiv a\pmod N$, there is a unital $K$-linear $*$-isomorphism $\Phi:R\longrightarrow M_d(R)$ such that $\Phi(e_i^*)\in M_d(R_{\le2})$ for each $i$.
\end{proposition}

\begin{proof}
Write $n=Qa+r$ and $\beta=\beta_a=w+1$. In the path on $S_c$, let $p_{c,i}$ be the vertex at distance $i$ from its start, so $p_{1,0}=1$ and $p_{2,0}=a$. For an integer $M\ge2(a-1)$, set
$$
d_{p_{c,i}}=M+c-1-2i\ge0
\qquad\text{and}\qquad
D_{c,i}=4-c+2i.
$$
By Lemma~\ref{lem:forest}, the nonzero expanded entries arising from row $p_{c,i}$ have degrees
$$
\begin{array}{l|c}
\text{entry and column}&\text{degree after expansion}\\ \hline
v_0^*\text{ in column }p_{c,i+1}&d_{p_{c,i}}-d_{p_{c,i+1}}=2\\
e_u^*\text{ in column }1&1+d_{p_{c,i}}-d_1=c-2i\\
\lambda^*\text{ in column }a&|\lambda|+d_{p_{c,i}}-d_a=|\lambda|+2-D_{c,i}
\end{array}.
$$
The first row of the table applies when the path has a next vertex. Thus it suffices to assign leaves with $|\lambda|\le D_{c,i}$ to row $p_{c,i}$. If $r=a$, then $\delta=m=a-1$ and $S_2=\{a\}$, with $D_{2,0}=2$.

Prescribe the numbers $I_j$ of internal paths of length $j$ by a finitely supported sequence of nonnegative integers with
$$
I_0=1,\qquad I_j\le\beta I_{j-1}\quad(j\ge1),\qquad \sum_{j\ge0}I_j=a-1.
$$
At depth $j$, attach all children of the $I_{j-1}$ internal paths $\nu$ and choose $I_j$ internal vertices among their $\beta I_{j-1}$ extensions $\nu e_u$ with $u\in B_a$. This produces a compatible tree. Let $\Lambda_c$ be its leaves of color $c$ and $F_c(D)=|\{\lambda\in\Lambda_c:|\lambda|\le D\}|$. For $D\ge1$,
$$
\begin{aligned}
F_1(D)&=\sum_{j=1}^D(\beta I_{j-1}-I_j)
=(\beta-1)\sum_{j=0}^{D-1}I_j+1-I_D,\\
F_2(D)&=(n-\beta)\sum_{j=0}^{D-1}I_j.
\end{aligned}
$$
Order the leaves $(\lambda_{c,k})_k$ of each color by nondecreasing length. Put $P_c(-1)=0$ and $P_c(i)=\sum_{q=0}^i C_a(p_{c,q})$. Assign $\lambda_{c,k}$ to row $p_{c,i}$ when $P_c(i-1)<k\le P_c(i)$. By Lemma~\ref{lem:counts}, this uses all leaves and positions. The length bounds hold if
$$
F_c(D_{c,i})\ge P_c(i)\qquad\text{for every row }p_{c,i}.
$$
This is automatic when $D_{c,i}$ is at least the tree height, since $F_c(D_{c,i})=|\Lambda_c|\ge P_c(i)$.

If $a\le w^2+3w+4=\beta^2+\beta+2$, put $h=a-2$ and $t=\min(\beta,h)$. Take $(I_0,I_1,I_2)=(1,t,h-t)$ and $I_j=0$ for $j\ge3$. Then $I_1\le\beta$, $\sum_{j\ge0}I_j=a-1$, and
$$
\begin{aligned}
h\le\beta&\ \Longrightarrow\ I_2=0,\\
h>\beta&\ \Longrightarrow\ I_2=h-\beta\le\beta^2=\beta I_1.
\end{aligned}
$$
All leaf lengths are at most three. Only row $a=p_{2,0}$ has $D_{c,i}<3$. If $I_2=0$, its bound $D_{2,0}=2$ is automatic. Otherwise $I_1=\beta$ and $\beta\le n-1$ by Lemma~\ref{lem:partition}, so
$$
F_2(2)=(n-\beta)(1+\beta)=n+\beta(n-\beta-1)\ge n>C_a(a)=P_2(0).
$$

If $a>w^2+3w+4$ and $N\le w^3$, then $w\ge3$ since $a<N$ and $w^3\le w^2+3w+4$ for $w=1,2$. Set
$$
E=a-\beta^2-\beta-2>0,\qquad z=\Bigl\lceil\frac E{\beta+1}\Bigr\rceil,\qquad t=E-z.
$$
Take $(I_0,I_1,I_2,I_3,I_4)=(1,\beta,\beta^2,z,t)$ and $I_j=0$ for $j\ge5$. The constraints on $(I_j)$ follow from
$$
\begin{aligned}
1\le z&\le E\le(\beta+1)z,\qquad 0\le t\le\beta z,\\
z&\le\frac{N}{w+2}+1\le\frac{w^3}{w+2}+1\le w^2<\beta^3,\\
\sum_{j\ge0}I_j&=1+\beta+\beta^2+E=a-1.
\end{aligned}
$$
All leaf lengths are at most five. Since $f(f(a))=a\ge2$ and $f(N-1)=1$, we have $w\le N-2$. As $C_a(p)\le N$ for each row $p$, the rows with $D_{c,i}<5$ satisfy
$$
\begin{aligned}
F_1(3)&=\beta^3-z\ge w^3+2w^2+3w+1>N\ge P_1(0),\\
F_2(2)-2N&=(N-w)(w+2)-2N=w(N-w-2)\ge0,\\
P_2(i)&\le(i+1)N\le F_2(2)\le F_2(D_{2,i})
\qquad\text{for }0\le i<\min(2,|S_2|).
\end{aligned}
$$

In either case, Proposition~\ref{prop:canonical} supplies a Leavitt generating family $(A_1,\ldots,A_n)$ in $M_a(R)$. Each row of each $A_i$ has a single nonzero entry, a path adjoint. Let $F$ be the expansion from Lemma~\ref{lem:forest} for the chosen depths. Then
$$
B_i=F(A_i)\in M_e(R_{\le2}),\qquad e=\sum_{p=1}^a n^{d_p}\equiv a\pmod N.
$$
The row property in Lemma~\ref{lem:forest} ensures the support condition, so Lemma~\ref{lem:lifting} supplies a Leavitt generating family $(X_1,\ldots,X_n)$ with $X_i\in M_d(R_{\le2})$ for every positive $d\equiv a\pmod N$. The assignments $e_i^*\mapsto X_i$ and $e_i\mapsto X_i^*$ induce the asserted isomorphism, with $X_i^*\in M_d(R_{\ge-2})$.
\end{proof}

To reach further residue classes, apply the inner isomorphism $\psi:R\to M_b(R)$ entrywise to the outer matrices of size $a$ from Proposition~\ref{prop:canonical}. The inner tree has height $\le2$, and we permute the edge indices. Each nonzero entry within a block is a path adjoint of length $\le4$. For a leaf block $\psi(\lambda^*)$ in outer row $p$, the bound in inner row $z$ is $4-\mathbf1_{\{1,b\}}(z)$ if $p=1$ and $2$ if $p=a$. By Lemma~\ref{lem:forest}, we choose depths so that the expanded ghost-edge images have entries in $R_{\le2}$.

\begin{proposition}\label{prop:composition}
Let $3\le a,b<n$ satisfy $\gcd(a,N)=\gcd(b,N)=1$. Write $n=Q_hh+r_h$ for $h\in\{a,b\}$, and assume $2\le r_h<h$ for $h\in\{a,b\}$, $Q_a\le Q_b$, and $\beta_b\ge b-2$.

Put $\beta=\beta_a$, $H=a-2$, $C=C_a(1)$, and $\Delta=\beta-Q_a-1$. Put $L_1=\beta_b-b+2$ and $L_2=n-\beta_b$. Set $K_0=Q_b-Q_a+1+\min(L_1,L_2)$, $u=\min(\beta,K_0)$, and $k=\min(u,H)$.

Assume the following three inequalities hold.
\begin{enumerate}[leftmargin=*,nosep]
\item $H\le 2\beta+u(u-1)$.
\item $\beta+H(u-1)+(\beta-u)k\ge C$.
\item $0\le\Delta\le n-\beta$.
\end{enumerate}
Then for every field $K$ and every positive integer $d\equiv ab\pmod N$, there is a unital $K$-linear $*$-isomorphism $\Phi:R\longrightarrow M_d(R)$ such that $\Phi(e_i^*)\in M_d(R_{\le2})$ for each $i$.
\end{proposition}

\begin{proof}
For size $b$, take $v_0$ and $b-2$ distinct paths $e_j$ with $j\in B_b$ as the internal vertices of a compatible tree. It has height at most two, with $L_1=\beta_b-(b-2)$ and $L_2=n-\beta_b$ leaves of length one of colors one and two, respectively.

By Lemma~\ref{lem:partition}, $Q_h+1\le\beta_h\le n-Q_h-1$ for $h\in\{a,b\}$. Since $K_0\ge1$, we have $1\le u\le\beta$. Choose $U\subseteq B_a$ with $|U|=u$ and fix $v\in U$. Since $r_a<a$, the only term of $X_{Q_a+1}$ in row $a$ is $E_{a,a-r_a+1}$, with $a-r_a+1\ge2$. By Proposition~\ref{prop:canonical}, the ghost edges in outer position $(a,1)$ are therefore indexed by $D=\{a,2a,\ldots,Q_aa\}\subseteq B_a^c$.

Put $J=\{Q_b+3,\ldots,n\}$. Map $D$ bijectively onto $\{1,\ldots,Q_a\}$ and send $v$ to $Q_b+1$. Assign distinct indices to $U\setminus\{v\}$, using $\{Q_a+1,\ldots,Q_b\}$ first and then $J$. The capacity bounds are
$$
\begin{aligned}
\max\{0,u-1-(Q_b-Q_a)\}&\le\min\{L_1,L_2\},\\
L_1&\le n-Q_b-b+1\le |J|\le C_b(1),\\
L_2&\le n-Q_b-1=C_b(b).
\end{aligned}
$$
Since $D\cap U=\varnothing$, these assignments extend to a permutation $\pi$. For each $j\in\pi(U)\cap J$, reserve a leaf of length one of each color for positions $(j,1)$ and $(j,b)$, respectively, using distinct leaves within each color. Assign the remaining leaves of length one of color two to unused positions in row $b$. Complete the bijection by color using Lemma~\ref{lem:counts}. Let $\phi_b$ be the isomorphism of Proposition~\ref{prop:canonical} for this assignment and set $\psi=\phi_b\circ\alpha_\pi$, where $\alpha_\pi$ is the $*$-automorphism $e_j\mapsto e_{\pi(j)}$.

For a matrix $W$ whose sole nonzero entry in each row is a path adjoint, let $t_z(W)$ and $\ell_z(W)$ be the column and path length of that entry in row $z$, with $1=v_0^*$ of length zero. The product rule following Proposition~\ref{prop:canonical} yields
$$
t_z(AB)=t_{t_z(A)}(B)\quad\text{and}\quad
\ell_z(AB)=\ell_z(A)+\ell_{t_z(A)}(B).
$$
Put $E_0=\{1,b\}$ and $\mu_z=\mathbf1_{E_0}(z)$. Let $G$ range over the matrices $\psi(e_j^*)$, let $S=\psi(e_v^*)$, and let $V$ range over the matrices $\psi(e_j^*)$ with $j\in U\setminus\{v\}$. Write $P$ for $S$ or any such $V$, with subscripts indicating independent choices. By the formulas of Proposition~\ref{prop:canonical} and the leaf assignments,
$$
\begin{aligned}
\ell_z(G)>0&\ \Longrightarrow\ (t_z(G),\ell_z(G))\in\{(1,1),(b,1),(b,2)\},\\
\ell_z(S)&=\mathbf1_{\{1,\ldots,r_b\}}(z),\\
1&\le\ell_z(V)\le2-\mu_z,\qquad \ell_z(P)\le2-\mu_z.
\end{aligned}
$$
For $V=\psi(e_j^*)$, the bound follows from the column-$1$ ghost edges if $\pi(j)\le Q_b$ and the length-one reservations in rows $1$ and $b$ if $\pi(j)\in J$.

Use (a) to construct the outer tree and (b), (c) to assign leaves to rows $1$ and $a$, respectively.
\begin{enumerate}[leftmargin=*,nosep]
\item Put $t=\min(H,\beta)$. Since $k=\min(u,t)\ge1$ and $t-k\le\beta-u$, choose $T\subseteq B_a$ with $v\in T$, $|T|=t$, and $|T\cap U|=k$. Take $e_i$ with $i\in T$ as the internal paths of length one. The candidates for internal paths of length two are
$$
\{e_v e_j:j\in B_a\}\sqcup
\{e_i e_j:i\in(T\cap U)\setminus\{v\},\ j\in U\}.
$$
The set has $\beta+(k-1)u$ elements. If $H\le\beta$, then $H-t=0$. Otherwise $(t,k)=(\beta,u)$ and (a) implies $H-t\le\beta+u(u-1)=\beta+(k-1)u$. Choose $H-t$ of these paths for $\mathcal I$. Together with $v_0$ and the $e_i$ for $i\in T$, they are the $a-1$ internal vertices of a compatible tree of height at most three.

\item Before expanding $\mathcal I$, select the color-one leaves
$$
\begin{aligned}
\Lambda_0
&=\{e_j:j\in B_a\setminus T\}\sqcup\{e_i e_j:i\in T,\ j\in U\}\\
&\quad{}\sqcup\{e_i e_j:i\in T\cap U,\ j\in B_a\setminus U\}.
\end{aligned}
$$
Since $\mathcal I\subseteq\Lambda_0$, the selection after expansion is $\Lambda_1=(\Lambda_0\setminus\mathcal I)\sqcup\{\nu e_j:\nu\in\mathcal I,\ j\in U\}$. By (b),
$$
\begin{aligned}
|\Lambda_1|&=\beta-t+ut+k(\beta-u)+(u-1)(H-t)\\
&=\beta+H(u-1)+(\beta-u)k\ge C.
\end{aligned}
$$
Assign $C$ of these leaves to row $1$.

\item By (c), choose $J_a\subseteq B_a^c$ with $|J_a|=\Delta$. The paths in
$$
\{e_j:j\in B_a^c\}\sqcup\{e_v e_j:j\in J_a\}
$$
are leaves of color two, since $e_v$ is internal and each nonroot internal path ends in an edge indexed by $B_a$. Assign these $(n-\beta)+\Delta=C_a(a)$ leaves to row $a$.
\end{enumerate}
Complete the bijection by color using Lemma~\ref{lem:counts}.

For $\lambda=e_{j_1}\cdots e_{j_s}$ put $W_\lambda=\psi(\lambda^*)=\psi(e_{j_s}^*)\cdots\psi(e_{j_1}^*)$. The leaves of length two or three selected in (b) yield
$$
\begin{array}{c|l|c}
\lambda&\text{indices}&W_\lambda\text{ of the form}\\ \hline
e_i e_j&i\in T,\ j\in U&PG\\
e_i e_j&i\in T\cap U,\ j\in B_a\setminus U&GP\\
e_v e_j e_\ell&j\in B_a,\ \ell\in U&PGS\\
e_i e_j e_\ell&i\in(T\cap U)\setminus\{v\},\ j,\ell\in U&P_1P_2V
\end{array}.
$$
In (c), $W_\lambda$ has the form $G$ or $GS$. A leaf of length three from (a) yields $G_1G_2S$ if it begins with $e_v$ and $GPV$ otherwise.

If the first factor has positive length in row $z$, its column lies in $E_0$, where $\ell_y(V)=1$ and $\ell_y(P)\le1$. For $GS$, the inequality $2\le r_b<b$ implies $(\ell_1(S),\ell_b(S))=(1,0)$, and $\ell_z(G)=2$ forces $t_z(G)=b$. The product rule yields
$$
\begin{aligned}
\ell_z(GS)&\le\max\{1,1+\ell_1(S),2+\ell_b(S)\}=2,\\
\ell_z(PV)&\le\max\{2,\ell_z(P)+1\}\le3-\mu_z,\\
\ell_z(GP)&\le\max\{2,\ell_z(G)+1\}\le3.
\end{aligned}
$$
In particular, $\ell_y(PV)\le2$ for $y\in E_0$. The same product rule yields
$$
\begin{aligned}
\ell_z(GPV)&\le\max\{3,\ell_z(G)+2\}\le4,\\
\ell_z(P_1P_2V)&\le\max\{3,\ell_z(P_1)+2\}\le4-\mu_z,\\
\ell_z(PG),\ \ell_z(PGS)&\le(2-\mu_z)+2,\\
\ell_z(G_1G_2S)&\le2+2.
\end{aligned}
$$
For $|\lambda|\le2$, the general bound is $\ell_z(W_\lambda)\le2|\lambda|\le4$. Thus for a leaf assigned to outer row $p$,
$$
\ell_z(W_\lambda)\le
\begin{cases}
4-\mu_z,&p=1,\\
2,&p=a,\\
4,&p\notin\{1,a\}.
\end{cases}
$$

Let $\phi_a$ be the isomorphism of Proposition~\ref{prop:canonical} for the outer assignment. With $M_a(M_b(R))\cong M_{ab}(R)$ indexed lexicographically, set
$$
\Theta=\operatorname{Mat}_a(\psi)\circ\phi_a:R\longrightarrow M_{ab}(R).
$$
Here $\operatorname{Mat}_a(\psi)$ applies $\psi$ entrywise. The substitution rule following Proposition~\ref{prop:canonical} ensures that each row of each $\Theta(e_j^*)$ has a unique nonzero entry, a path adjoint.

On the outer row-index paths $(p_0,\ldots,p_{m-1})$ and $(q_0,\ldots,q_{a-m-1})$, with $p_0=1$ and $q_0=a$, take an integer $M\ge2(a-1)$ and set
$$
d_{(p_i,z)}=M+\mu_z-2i\quad\text{and}\quad
d_{(q_i,z)}=M+2-2i.
$$
These depths are nonnegative and $d_{(p,z)}\le M$ for $p\notin\{1,a\}$. By Lemma~\ref{lem:forest}, a nonzero expanded entry arising from a path adjoint of length $s$ at $((p,z),(q,z'))$ has degree $s+d_{(p,z)}-d_{(q,z')}$. A scalar block $I_b$ preserves $z$, so along the directed paths this degree is $d_{(p,z)}-d_{(q,z)}=2$. For leaf blocks the length bounds imply
$$
\ell_z(W_\lambda)+d_{(p,z)}\le M+4=d_{(a,t_z(W_\lambda))}+2.
$$
For $G=\psi(e_j^*)$ in outer column $1$, put $z'=t_z(G)$. Since positive lengths end in $E_0$, $\ell_z(G)-\mu_{z'}\le1$. If $p=a$, then $j\in D$ and $\pi(j)\le Q_a\le Q_b$, so $(\ell_z(G),z')=(1,1)$. Hence
$$
\ell_z(G)+d_{(p,z)}-d_{(1,z')}\le
\begin{cases}
1+\mu_z\le2,&p=1,\\
1+(M+2)-(M+1)=2,&p=a,\\
2+M-M=2,&p\notin\{1,a\}.
\end{cases}
$$

The expansion $F$ of Lemma~\ref{lem:forest} at these depths yields a Leavitt generating family $(B_i)_{i=1}^n$ with
$$
B_i=F(\Theta(e_i^*))\in M_e(R_{\le2}),\qquad
e=\sum_{p,z}n^{d_{(p,z)}}\equiv ab\pmod N.
$$
The row property in Lemma~\ref{lem:forest} ensures the support condition in Lemma~\ref{lem:lifting}. For every $d\ge1$ with $d\equiv ab\pmod N$, the latter lemma yields a Leavitt generating family $(X_i)_{i=1}^n$ in $M_d(R)$ with $X_i\in M_d(R_{\le2})$ and $X_i^*\in M_d(R_{\ge-2})$. The assignments $e_i^*\mapsto X_i$ and $e_i\mapsto X_i^*$ induce $\Phi$.
\end{proof}

\begin{remark}\label{rem:boundary-eleven}
For $n=11$ and $a=b=3$, we have $Q_a=Q_b=3$, $r_a=r_b=2$, $\beta_a=\beta_b=4$, $H=1$, $C=7$, $\Delta=0$, $L_1=3$, $L_2=7$, $K_0=u=4$, and $k=1$. The three inequalities of Proposition~\ref{prop:composition} reduce to
$$
1\le8+12,\qquad 4+1(4-1)=7,\qquad 0\le0\le7.
$$
Thus the proposition applies to $ab=9\equiv-1\pmod{10}$. In each tree, the internal vertices are $v_0$ and one path of length one.
\end{remark}

\section{Proof of the main theorem}\label{sec:remaining}\label{sec:powers}

For the standard Leavitt generating family $(e_i^*)_{i=1}^n$ in $M_1(R)$, we have $\mathcal S_1=\varnothing$ since $e_i^*\in R_1$ for all $i$. Lemma~\ref{lem:lifting} transfers the bound to $d\ge1$ with $d\equiv1\pmod N$, including every $d\ge1$ when $n=2$.

\paragraph{Prime $n$.} For a unit $w$ modulo $N$ with $w\ge3$ and $w^3<N$, choose $b=w^j$ and $a=f(wb)$ so that the sizes satisfy Proposition~\ref{prop:composition} after ordering, with $ab\equiv f(w)\pmod N$.

\begin{proposition}\label{prop:powers}
Suppose $n$ is prime, and let $w$ be a unit modulo $N$, represented in $\{1,\ldots,N-1\}$, with $w\ge3$ and $N>w^3$. Choose the least integer $j\ge1$ for which $w^j\ge\sqrt N/w$, and put $b=w^j$, $A=wb$, and $a=f(A)$.

The ordered pair $(\max\{a,b\},\min\{a,b\})$ satisfies the hypotheses of Proposition~\ref{prop:composition}, with its first and second coordinates in the roles of $a$ and $b$, respectively. Moreover, $ab\equiv f(w)\pmod N$. Consequently, every positive dimension $d\equiv f(w)\pmod N$ admits a Leavitt generating family $(X_1,\ldots,X_n)$ with $X_i\in M_d(R_{\le2})$.
\end{proposition}

\begin{proof}
Put $T=\sqrt N$. Since $T/w>1$ and $j$ is minimal, we have $w^{j-1}<T/w\le b<T$ and $T\le A$. If $j=1$, then $wA=w^3<N$. If $j\ge2$, then $wA=w^{j+2}\le w^{2j}<N$. Thus
$$
T\le A<N/w\le N/3.
$$
Both $b$ and $A$ are representatives of units. Since $A\ne N-1$, we have $a\ne1$. Since $N\ge28$, the prime $n$ is odd and $N$ is even, so the unit $a$ is not two. Hence $3\le a,b<n$. By the convention of Section~\ref{sec:partition}, $2\le r_h\le h$ for $h\in\{a,b\}$. Since $1<h<n$ and $n$ is prime, $h\nmid n$, so $r_h<h$. As $aA\equiv-1\pmod N$ and $A=wb$, we have
$$
f(a)=A,\qquad ab\equiv f(w)\pmod N,\qquad B:=f(b)=[wa]_N,
$$
where $[wa]_N$ is the representative in $\{1,\ldots,N-1\}$. Put $D=N-B$. Then $bB=m_bN-1$ and $bD=t_bN+1$ for integers $m_b,t_b\ge1$, since $b\ge3$ and $B,D\ge1$. As $b^2\le N-1$,
$$
B\ge (N-1)/b\ge b\quad\text{and}\quad D\ge (N+1)/b>T.
$$
When $b$ is the inner dimension, $\beta_b=B+1\ge b-2$ ensures that a compatible tree of height at most two exists. The numbers of leaves of length one of colors one and two are $L_1=B-b+3\ge3$ and $L_2=D>T$, respectively.

Write $\mathcal R=\beta+H(u-1)+(\beta-u)\min(u,H)$ in the notation of Proposition~\ref{prop:composition}. If $u\ge H$ or $u=\beta$, then $\mathcal R=\beta+H(\beta-1)$. By Lemma~\ref{lem:counts}, this is the total number of leaves of color one, which is at least the number of monomial positions in row $1$ of the outer construction.

Suppose $a\le b$, so $Q_b\le Q_a$. The outer size is $b$ and the inner size is $a$. Since $\beta_a=A+1\ge a-2$, a compatible inner tree of height at most two exists. Its numbers of leaves of length one satisfy
$$
A-a+3\ge(w-1)b+3>b-2\quad\text{and}\quad N-A>2N/3>T>b-2.
$$
Consequently $K_0,\beta_b>H=b-2$, so $u\ge H$, which implies $H\le 2\beta_b+u(u-1)$ and $\mathcal R\ge C_b(1)$. Since $wa\le A<N/3$, we have $B=wa$ and $2B<N$. Using $B\ge Q_b$ from Lemma~\ref{lem:partition}, we obtain $0\le B-Q_b\le N-B$, the last-row inequality for $(b,a)$.

Suppose $a>b$, so $Q_a\le Q_b$. With outer size $a$ and inner size $b$, we have $\beta=A+1$, $H=a-2$ and $u=\min(\beta,K_0)$, where $K_0=Q_b-Q_a+1+\min(L_1,L_2)$. By Lemma~\ref{lem:partition}, $A\ge Q_a$. Since $2A<N$, we obtain $0\le A-Q_a\le N-A$.

If $a>A+3$, then $a\ge3b+4$ and $H>\beta$. For $u=\beta$, we have
$$
H<a<N\le A^2<2\beta+u(u-1)\quad\text{and}\quad \mathcal R\ge C_a(1).
$$
If $u=K_0<\beta$, then $u\ge Q_b-Q_a+4$ since $L_1\ge3$ and $D>T>3$. Since $Q_a a\le N-1$ and $N=Q_b b+r_b-1<(Q_b+1)b$, we obtain
$$
3Q_a b<Q_a a\le N-1<(Q_b+1)b\quad\text{and}\quad b^2<N<(Q_b+1)b.
$$
Thus $3Q_a\le Q_b$, $b\le Q_b$, $3u\ge2Q_b+12$ and $N<(Q_b+1)^2$. If $Q_a\le2$, then $u^2\ge(Q_b+2)^2>N>a$. If $Q_a\ge3$, then $3a\le N-1$ and
$$
9u^2\ge(2Q_b+12)^2>4(Q_b+1)^2>4N>9a.
$$
Hence $2\beta+u(u-1)>u^2>a>H$. Since $H\ge3b$ and $3(u-1)\ge2Q_b+9$,
$$
\mathcal R\ge H(u-1)\ge b(2Q_b+9)>(Q_b+1)b>N\ge C_a(1).
$$

It remains to consider $b<a\le A+3$. If $wa\ge N$, the identity $aA+1=m_aN$ from Lemma~\ref{lem:partition} implies $m_a\ge b+1$ because $aA=bwa$. Since $N\ge wA+1$ and $wA+b>3A$,
$$
m_aN\ge(b+1)(wA+1)=A^2+wA+b+1>A^2+3A+1\ge aA+1,
$$
a contradiction. Thus $wa<N$.

Since $B=wa$, we have $L_1=wa-b+3>(w-1)a+3>H$. Also $H\le\beta$, so $H\le2\beta+u(u-1)$. If $K_0\ge H$, then $u\ge H$ and $\mathcal R\ge C_a(1)$.

If $K_0<H$, then $\min(L_1,D)=D$ because $L_1>H$ and $Q_b-Q_a+1\ge1$. Hence $D+1\le u=K_0<H\le A+1$. Set $P(v)=\beta-H+v(H+\beta-v)$. We have $\mathcal R=P(u)$ and $P(v+1)-P(v)=H+\beta-2v-1\ge1$ for $1\le v<H$. Thus
$$
\mathcal R\ge P(D+1)=aD+(D+2)A-D^2-3D+1\ge N+4>C_a(1),
$$
since $aD\ge bD\ge N+1$ and $(D+2)A-D^2-3D+1\ge3$ by $A\ge D+1$.

Proposition~\ref{prop:composition} now applies to the stated ordered pair and yields the required family over every field and in every positive dimension congruent to $ab$ modulo $N$.
\end{proof}

For example, take $n=71$ and $d=23$, so $N=70$ and $w=f(23)=3$. The sufficient conditions in Proposition~\ref{prop:direct} fail, since $23>3^2+3\cdot3+4=22$ and $70>3^3$. Proposition~\ref{prop:powers} selects $j=1$, $b=3$, $A=9$ and $a=f(9)=31$, with $ab=93\equiv23\pmod{70}$. In the notation of Proposition~\ref{prop:composition},
$$
\begin{aligned}
(Q_a,r_a,\beta_a)&=(2,9,10),& (Q_b,r_b,\beta_b)&=(23,2,24),\\
(H,C,\Delta)&=(29,67,7),& (L_1,L_2,K_0)&=(23,47,45),
\end{aligned}
$$
so $u=k=10$. The three inequalities in Proposition~\ref{prop:composition} become
$$
29\le110,\qquad 271\ge67\quad\text{and}\quad 0\le7\le61.
$$
The outer tree has $t=10$ internal vertices of depth one and $H-t=19$ of depth two. The construction yields a Leavitt generating family in $M_{23}(L_K(1,71))$ with the required degree bound.

For the class $-1\bmod N$, we choose explicit composition sizes when $n\ge11$ is prime.

\begin{proposition}\label{prop:minusone}
If $n$ is prime, every positive $d\equiv-1\pmod N$ admits a Leavitt generating family $(X_1,\ldots,X_n)$ with $X_i\in M_d(R_{\le2})$.
\end{proposition}

\begin{proof}
For $n=2,3$, we have $-1\equiv1\pmod N$. For $n=5,7$, take $a=n-2$, so $f(a)=1$ and $a\le8$, and apply Proposition~\ref{prop:direct}. For $n=11$, take $a=b=3$ as in Remark~\ref{rem:boundary-eleven}.

Use the notation of Proposition~\ref{prop:composition}, with $\mathcal R=\beta+H(u-1)+(\beta-u)\min(u,H)$. For the remaining primes, put
$$
(a,b)=
\begin{cases}
(2\kappa+1,2\kappa-1),&n=4\kappa+1,\quad \kappa\ge3,\\
(4\kappa-1,2\kappa+1),&n=8\kappa+3,\quad \kappa\ge2,\\
(6\kappa+5,4\kappa+1),&n=8\kappa+7,\quad \kappa\ge2.
\end{cases}
$$
Here $\kappa$ is an integer. In the first two cases $ab+1=\kappa N$, and in the third $ab+1=(3\kappa+1)N$. Thus $a,b$ are units modulo $N$, with $f(a)=b$ and $f(b)=a$. The parameters are
$$
\begin{array}{c|cccc|ccc}
n&Q_a&r_a&Q_b&r_b&\beta&H&K_0\\ \hline
4\kappa+1&1&2\kappa&2&3&2\kappa&2\kappa-1&7\\
8\kappa+3&2&5&3&2\kappa&2\kappa+2&4\kappa-3&2\kappa+3\\
8\kappa+7&1&2\kappa+2&2&5&4\kappa+2&6\kappa+3&2\kappa+3
\end{array}.
$$
It follows that $3\le b\le a<n$, $2\le r_a<a$, $2\le r_b<b$ and $Q_a\le Q_b$. Also $\beta_b=a+1\ge b-2$ by Lemma~\ref{lem:partition}, and $H<2\beta\le2\beta+u(u-1)$. The last-row inequalities are
$$
\begin{aligned}
0&\le2\kappa-2\le2\kappa+1&&\text{if }n=4\kappa+1,\\
0&\le2\kappa-1\le6\kappa+1&&\text{if }n=8\kappa+3,\\
0&\le4\kappa\le4\kappa+5&&\text{if }n=8\kappa+7,
\end{aligned}
$$
where the middle and right quantities are $\Delta$ and $n-\beta$, respectively.

If $n=4\kappa+1$ with $\kappa=3,4$, then $u=\min(2\kappa,7)\ge H$. If $n=8\kappa+3$, then $u=\beta$. In both cases $\mathcal R=\beta+H(\beta-1)\ge C_a(1)$ by Lemma~\ref{lem:counts}. For $n=4\kappa+1$ with $\kappa\ge5$ and for $n=8\kappa+7$, respectively,
$$
\begin{aligned}
\mathcal R&=2\kappa+6(2\kappa-1)+7(2\kappa-7)\\
&=28\kappa-55\ge4\kappa-2=C_a(1),\\
\mathcal R&=(4\kappa+2)+(6\kappa+3)(2\kappa+2)+(2\kappa-1)(2\kappa+3)\\
&=16\kappa^2+26\kappa+5\ge8\kappa+4=C_a(1).
\end{aligned}
$$
Among odd integers $n\ge3$, the three ranges omit only $3,5,7,9,11,15$. Their prime members were treated above. In each remaining case, Proposition~\ref{prop:composition} applies, with $ab\equiv-1\pmod N$.
\end{proof}

\begin{example}
For $n=43=8\cdot5+3$, the pair is $(a,b)=(19,11)$, with $ab=209=5\cdot42-1$. Proposition~\ref{prop:composition} therefore supplies a generating family in $M_{41}(L_K(1,43))$ and in every positive dimension congruent to $41$ modulo $42$.
\end{example}

\paragraph{Composite $n$.} If $p\ge2$ and $p\mid n$, it follows from \cite[Prop.~6.3]{AAP} that $R\cong M_p(R)$ as graded algebras. Iterating it entrywise yields graded isomorphisms $M_d(R)\cong M_{dp^\ell}(R)$ for $\ell\ge0$.

The row matrices that Paschke and Salinas construct in \cite[Prop.~2.5]{PS1979} become $Z_{c,j}^*$ when the Cuntz generators are replaced by $e_1,\ldots,e_n$. Their adjoints $Z_{c,j}$ are the images of the ghost edges.

\begin{lemma}\label{lem:divisor}
If $p\ge2$ divides $n$, then there is a graded unital $K$-linear $*$-isomorphism $\Gamma_p:R\longrightarrow M_p(R)$. For every $d\ge1$ and $\ell\ge0$, its successive entrywise applications, under the natural identifications $M_q(M_p(R))\cong M_{qp}(R)$, induce a graded unital $K$-linear $*$-isomorphism $M_d(R)\longrightarrow M_{dp^\ell}(R)$. The inverse preserves upper degree bounds and Leavitt generating families.
\end{lemma}

\begin{proof}
Write $n=pL$ and let $u_j(r)=(j-1)p+r$ for $1\le j\le L$ and $1\le r\le p$. For $1\le c\le p$, put
$$
Z_{c,j}=\sum_{r=1}^{p}e_{u_j(r)}^*E_{rc}
\qquad\text{and}\qquad
Z_{c,j}^*=\sum_{r=1}^{p}e_{u_j(r)}E_{cr}.
$$
By (CK1) and (CK2),
$$
\begin{aligned}
Z_{c,j}Z_{c',j'}^*
&=\delta_{cc'}\sum_{r,s=1}^p e_{u_j(r)}^*e_{u_{j'}(s)}E_{rs}
=\delta_{cc'}\delta_{jj'}I_p,\\
\sum_{j=1}^L Z_{c,j}^*Z_{c',j}
&=\sum_{j=1}^L\sum_{r=1}^p e_{u_j(r)}e_{u_j(r)}^*E_{cc'}=E_{cc'}.
\end{aligned}
$$
In particular, $\sum_{c,j}Z_{c,j}^*Z_{c,j}=\sum_cE_{cc}=I_p$. Thus $e_{u_j(c)}^*\mapsto Z_{c,j}$ and $e_{u_j(c)}\mapsto Z_{c,j}^*$ define a unital $K$-linear $*$-homomorphism $\Gamma_p:R\to M_p(R)$.

Let $A=\operatorname{im}\Gamma_p$. The second identity shows that $E_{cc'}\in A$ for all $c,c'$. Hence
$$
E_{sr}Z_{1,j}E_{1t}=e_{u_j(r)}^*E_{st}\in A
\qquad\text{and}\qquad
(e_{u_j(r)}^*E_{st})^*=e_{u_j(r)}E_{ts}\in A.
$$
These elements and the matrix units generate $M_p(R)$. By simplicity, $\ker\Gamma_p=0$, as in Proposition~\ref{prop:canonical}. Since $Z_{c,j}\in M_p(R_1)$ and $Z_{c,j}^*\in M_p(R_{-1})$, $\Gamma_p$ preserves the degree of each path monomial $\alpha\beta^*$ and is graded.

Applying $\Gamma_p$ entrywise and identifying $M_q(M_p(R))$ with $M_{qp}(R)$ by lexicographic indexing defines a graded unital $K$-linear $*$-isomorphism $\operatorname{Mat}_q(\Gamma_p)$. With $\Theta_\ell:M_d(R)\to M_{dp^\ell}(R)$, the iteration is
$$
\Theta_0=\operatorname{id}_{M_d(R)}
\qquad\text{and}\qquad
\Theta_{\ell+1}=\operatorname{Mat}_{dp^\ell}(\Gamma_p)\circ\Theta_\ell.
$$
Write $\Gamma_p^{(\ell)}$ for this iteration with $d=1$. By lexicographic indexing, $\Theta_\ell=\operatorname{Mat}_d(\Gamma_p^{(\ell)})$. If $C=\sum_s C_s$ with $C_s\in M_d(R_s)$ and $D=\Theta_\ell(C)\in M_{dp^\ell}(R_t)$, then the degree-$s$ component $D_s$ satisfies
$$
0=D_s=\Theta_\ell(C_s)\ \Longrightarrow\ C_s=0
\qquad\text{for }s\ne t.
$$
Thus $\Theta_\ell^{-1}$ is graded and preserves upper degree bounds. As a unital $*$-isomorphism, it also preserves Leavitt generating families.
\end{proof}

Choose a prime $p\mid n$ and $\ell\ge0$ so that $dp^\ell\equiv c'\pmod N$, where $1\le c'<N$. We arrange that $c'=1$ or Proposition~\ref{prop:direct} applies to $a=c'$. Lemma~\ref{lem:divisor} then transfers the bound to $M_d(R)$.

\begin{proposition}\label{prop:composite}
If $n$ is composite, every positive $d$ relatively prime to $N$ admits a Leavitt generating family $(X_1,\ldots,X_n)$ with $X_i\in M_d(R_{\le2})$.
\end{proposition}

\begin{proof}
Choose a prime $p\mid n$ with $p\le\sqrt n$. Since $\gcd(p,N)=1$, $p$ is a unit modulo $N$. Let $c\in\{1,\ldots,N-1\}$ represent $d$ modulo $N$, put $t=f(c)$, and choose the least $k\ge0$ with $(p^kt)^2\ge N$. For $0\le j<k$, it follows from minimality that
$$
p^{j+1}t<p\sqrt N\le\sqrt{nN}<N+1.
$$
Since $p^{j+1}t$ is a unit modulo $N$, it cannot equal $N$. Together with $1\le t<N$, this puts every $p^jt$, $0\le j\le k$, in $\{1,\ldots,N-1\}$.

Put $t'=p^kt$ and $c'=f(t')$. Then $f(c')=t'$ and $c'\equiv cp^{-k}\pmod N$. If $c'=1$, use the standard family $(e_1^*,\ldots,e_n^*)$ and Lemma~\ref{lem:lifting}. If $c'\ge2$, the inequality $c'\le N-1<(t')^2+3t'+4$ allows us to apply Proposition~\ref{prop:direct}. Both cases yield families in every positive dimension congruent to $c'$ modulo $N$.

Choose $\ell\ge0$ with $p^\ell\equiv p^{-k}\pmod N$, using the finiteness of the unit group modulo $N$. Then $dp^\ell\equiv c'\pmod N$. By Lemma~\ref{lem:divisor}, the isomorphism
$$
\operatorname{Mat}_d(\Gamma_p^{(\ell)}):M_d(R)\longrightarrow M_{dp^\ell}(R)
$$
and its inverse preserve the involution and each homogeneous degree. Pulling the family back to $M_d(R)$ preserves the Leavitt relations and generation, with $X_i\in M_d(R_{\le2})$ and $X_i^*\in M_d(R_{\ge-2})$ for each $i$.
\end{proof}

\begin{proof}[Proof of Theorem~\ref{thm:main}]
If $d\equiv1\pmod N$, use $(e_1^*,\ldots,e_n^*)$ and Lemma~\ref{lem:lifting}. This includes all positive $d$ when $n=2$. If $n$ is composite, apply Proposition~\ref{prop:composite}.

For odd prime $n$ and $d\not\equiv1\pmod N$, let $c\in\{2,\ldots,N-1\}$ represent $d$ modulo $N$ and put $w=f(c)$. Since $N$ is even, $w$ is odd. We use
$$
\begin{cases}
\text{Proposition~\ref{prop:minusone}},& w=1,\quad c=f(1)=N-1,\\
\text{Proposition~\ref{prop:direct}},& w\ge3,\quad c\le w^2+3w+4\ \text{or}\ N\le w^3,\\
\text{Proposition~\ref{prop:powers}},& w\ge3,\quad c>w^2+3w+4,\quad N>w^3.
\end{cases}
$$
In the second case take $a=c$, and in the third use $f(w)=c$. Each case yields a Leavitt generating family $(X_1,\ldots,X_n)$ with $X_i\in M_d(R_{\le2})$. The adjoints $Y_i=X_i^*$ belong to $M_d(R_{\ge-2})$.
\end{proof}

\begin{remark}\label{rem:distinct-generators}
The $2n$ matrices in Theorem~\ref{thm:main} are distinct. If $i\ne j$, then $(X_i-X_j)Y_i=I_d= X_i(Y_i-Y_j)$, so the matrices are distinct within each family. If $X_i=Y_j$ with $i\ne j$, then $X_i^2=0$ and $X_jX_i=I_d$, whence $X_i=(X_jX_i)X_i=0$, contrary to $X_iY_i=I_d$. If $X_i=Y_i$, then $X_i^2=I_d$. For $j\ne i$, we then have $Y_j=X_i(X_iY_j)=0$, contrary to $X_jY_j=I_d$.
\end{remark}

\section*{Declaration of AI Use}
The authors used AI to support discussions of mathematical ideas during the preparation of this manuscript. They reviewed and revised the manuscript before finalizing it and take responsibility for its content.

\end{document}